\documentclass[letterpaper,10pt]{article}

\usepackage[letterpaper,left=1in,right=1in,top=1.5in,bottom=1.5in]{geometry}
\newcommand{\myauthor}{Jan Cao}
\newcommand{\mytitle}{}
\usepackage{tikz}
\usetikzlibrary{matrix,arrows}
\usetikzlibrary{babel}

\usepackage{setspace}
\usetikzlibrary{matrix,arrows}
\usepackage{tikz-cd}
\title{Mathematizing the Limits of Time:\\ Heidegger, Derrida, and the Topology of Temporality}
\author{Jan Cao}
 \date{\vspace{-5ex}}

\usepackage[pdfstartview=FitH,
            pdfauthor={\myauthor},
            pdftitle={\mytitle},
            colorlinks,
            linkcolor=reference,
            citecolor=citation,
            urlcolor=e-mail,
            backref]{hyperref}

\usepackage{amsmath}
\usepackage{amscd}
\usepackage{amsbsy}
\usepackage{amssymb}
\usepackage{verbatim}
\usepackage{eufrak}
\usepackage{eucal}
\usepackage{microtype}
\usepackage{hyperref}
\usepackage{mathrsfs}
\usepackage{amsthm}
\usepackage{stmaryrd}
\usepackage[all,cmtip]{xy}
\usepackage{tikz}
\usetikzlibrary{matrix,arrows}
\usepackage[abbrev,lite,alphabetic]{amsrefs}
\usepackage{dsfont}
\usetikzlibrary{babel}
\usepackage{fancyhdr}
\providecommand{\keywords}[1]{\textbf{\textit{Keywords---}} #1}
\usepackage{makeidx}
\makeindex

\usepackage{color}
\definecolor{todo}{rgb}{1,0,0}
\definecolor{conditional}{rgb}{0,1,0}
\definecolor{e-mail}{rgb}{0,.40,.80}
\definecolor{reference}{rgb}{.20,.60,.22}
\definecolor{mrnumber}{rgb}{.80,.40,0}
\definecolor{citation}{rgb}{0,.40,.80}

\let\oldmarginpar\marginpar
\renewcommand\marginpar[1]{\-\oldmarginpar[\raggedleft\footnotesize #1]%
{\raggedright\footnotesize #1}}

\usepackage{relsize}
\usepackage[bbgreekl]{mathbbol}
\usepackage{amsfonts}
\usepackage{dirtytalk}
\DeclareSymbolFontAlphabet{\mathbb}{AMSb} 
\DeclareSymbolFontAlphabet{\mathbbl}{bbold}

\begin{document}

\maketitle
\begin{abstract}
    
\say{The mathematization of time has limits,} writes Derrida in \say{Ousia and Gramme.} Taking this quote in all possible senses, this paper considers Derrida’s definition of limit as \emph{gramme}, trace, and aporia, and develops the \say{mathematization} of all three. I will consider the structure of mathematical limit in Fregean arithmetics, calculus, and topology. I argue that a topological approach to the concept of \emph{gramme} creates a new \say{limit} that illustrates the map of our \say{commute} to the unknown, and therefore avoids the aporetic act of crossing the impassible limit between life and death. This new \say{limit,} often referred to by mathematicians as a \say{cone,} is a shapeless shape that offers us knowledge about time and space that we otherwise cannot know about. By redefining the \say{limit} of time through three mathematical systems -- arithmetic, calculus, and topology -- this paper offers new perspectives to couple mathematics with philosophy. 

\end{abstract}

In \emph{The Basic Problems of Phenomenology}, a lecture course that Heidegger gave in 1927, he offers a critique of Aristotle’s concept of time, which states that time is the something \say{ counted} in connection with motion, something that shows itself in regards to the before and after(Aristotle, 1983, 219b1f.) Every \say{now} is in transition, and the transition is measured by number. However, Heidegger believes that each and every number does not reach the essential nature of time.(Heidegger 1988, p.249.) What is relevant in the Aristotelian concept of time is only its measurability, its countability; the essential nature of every singular now is, in a sense, irrelevant. But the problem is, Heidegger continues, how are the beings in time connected to the soul [Seele, \emph{psych{\=e}}]? In other words, how is time interconnected with the concept of the world and thus, with the structure of \emph{Dasein}? Heidegger's critique could be summarized as the following: the consequence of understanding time in the Aristotelian \say{vulgar} way is that, time becomes endless both in the past and in the future. It is also understood as \say{free-floatings \say{in-itself} of a course of \say{nows} which is present-at-hand}(Heidegger 1988, p.272). Since time becomes a succession of limitless nows, a ‘course of time’ that flows as an indifferent temporal medium, it is understood as merely a collection of objective facts. Every now becomes exactly the same as the next now. Therefore, time loses its datability [\emph{Datierbarkeit}], its worldhood [\emph{Weltlichkeit}], its spannedness [\emph{Gespanntheit}], its character of having a location [\emph{{\"O}rtlichkeit}], and so on(Heidegger 1988, 268). More importantly, time also loses its death-ability: it becomes infinite, and therefore cannot die its own death. This process is what Heidegger calls a \say{leveling off}[\emph{nivellieren}] or world time, in which our immediate experience as temporal beings becomes abstract, lifted off from its locus of ontological import(Derrida 1982, 35). This way of understanding time obscures the radical finitude of \emph{Dasein} that Heidegger emphasizes, as \emph{Dasein} must be understood as circumscribed by its own temporal ekstases. If there was endless time, \emph{Dasein} would lose its character of being-towards-death, and the meaninglessness of \emph{Dasein} will stretch into an indefinite future which will never be grasped authentically.

In \say{Ousia and Gramme,} an essay devoted to a single note in Heidegger's \emph{Being and Time}, Derrida points out that Heidegger's critique of Aristotle is still problematic, because the attempt to depart from the \say{vulgar concept of time} to an \say{authentic time} fails to escape the grammar and lexicon of metaphysics. The conceptual pairs of opposites that Heidegger discusses, such as being and non-being, limit and number, presence as \emph{Gegenw{\"a}rtigkeit} and presence as \emph{Anwesenheit}, are \say{ordered around one fundamental axis: that which separates the authentic from the inauthentic and, in the very last analysis, primordial from fallen temporality}(Derrida 1982, p.63). Derrida continues to question, \say{why qualify temporality as authentic -- or proper [\emph{eigentlich}]-- and as inauthentic -- or improper -- when every ethical preoccupation has been suspended?}(Derrida 1982, p.63) His answer is simple: we examine the concept of temporality in terms of these oppositions, only because our entire metaphysical structure is dependent on it. Heidegger’s thought attempts to solve the problem of time in the original Greco-Western-philosophy by assigning it a closure; however, the closure necessarily creates a new, uncrossable \say{limit} that undermines the very effort to mend a similar one.

\section{The Mathematization of Time ``Has Limits"}

According to Derrida's reading of Heidegger's critique, the Aristotelian mathematization of time \say{has limits}:
\begin{quote}
Time is not thought on the basis of the now. It is for this reason that the mathematization of time  \emph{has limits}... It is the extent to which time requires \emph{limits}, nows analogous to points, and in the extent to which the limits are always accidents and potentialities, that time cannot be made perfectly mathematical, that time's mathematization has limits, and remains, as concerns its essence, accidental(Derrida 1982, p.61).
\end{quote}
The mathematization of time -- that is, to think of time as counted \say{nows} -- has \say{limit} in three senses. The first exists in the Heideggerian sense, that the term \say{limit} indicates something like a gap between the \say{vulgar,} \say{mathematical} concept of time and primordial [\emph{urspr{\"u}ngliche}], authentic temporality. \say{Nows} as numbers do not belong to the essence of time; they are always elsewhere, accidental(Derrida 1982, p.61). Secondly, the mathematization of time has \say{limit} is Derrida's way of saying that, to think of time mathematically is \say{impossible,} because time is always elsewhere, and that the attempt to de-limit the \say{vulgar} concept of time is again limited by its own metaphysical, ontological nature. And yet time \say{requires limits}, or \say{nows,} that serve to measure and to enumerate(Derrida 1982, p.61).Years later, in \emph{Aporias}, Derrida again explicates this impossibility by posing the following questions: \say{What if there was no other concept of time than the one that Heidegger calls \say{vulgar}? What if, consequently, opposing another concept to the \say{vulgar} concept were itself impracticable, nonviable, and impossible?}(Derrida 1993, p.14) The fact that the mathematization of time is simultaneously required and impossible creates a \say{vulgar paradox}, which illuminates an \say{exoteric aporia} that both mathematics and metaphysics fail to cross(Derrida 1982, p.43).

Finally, \say{the mathematization of time} is only made possible through its own transgression, and therefore could only be thought of as a \say{trace} that is produced as its own erasure. Derrida's \say{trace} has no limit or boundary (\emph{p{\'e}ras}), as every \say{now} that constitutes it has ceased to be both beginning and end in its indefinite self-reproduction. It is a movement that goes in circles without reconstituting and re-enforcing predominant Western metaphysics and onto-theology. Instead, it functions as a critique of metaphysics from within: \say{only on this condition can metaphysics and our language signal in the direction of their own transgression}(Derrida 1982, p.61). The trace \say{traces} by following its own erasure. The trace is what makes it possible to understand time as simultaneously the now-no-longer and the now-not-yet, as both the \say{complete} \emph{gramme} that unfolds in space and the \say{incomplete} \emph{gramme} that moves within the teleology of the metaphysical system. \say{In this sense,} Derrida writes, \say{the circle removes the limit of the point only by developing its potentiality. The \emph{gramme} is \emph{comprehended} by metaphysics between the points and the circle, between potentiality and the act}(Derrida 1982, p.60).

According to Derrida, the mathematization of time also has \say{limit} in three different \say{shapes}: first, limit as a border that differentiates \say{vulgar time} from \say{authentic time,} or separates each \say{now} from the \say{nows} before and after it; second, limit as tracing (\emph{gramme, Linie}), a movement of transgression produced as its own erasure; and finally, limit as aporia, a limit beyond which there is nothing, at least nothing thinkable, to cross into(Derrida 1993, p.15). There exists an aporia between life and death that is impossible to be crossed, because we are unable to know what exists (or does not exist) on the other side of the aporia. On the one hand, aporia takes the \say{shape} of gift, forgiveness, or hospitality, all of which share the same paradoxical structure in which the very condition that makes it possible is simultaneously the same one that make it impossible. On the other hand, aporia is something unthinkable, something \say{shapeless,} precisely because of its paradoxical nature. Is impossible to cross or trespass aporia, since it is a \say{nonpassage because its elementary milieu does not allow for something that could be called passage, step, walk, gait, displacement, or replacement, a kinesis in general. There is no more path (\emph{odos, methodos, Weg}, or \emph{Holzweg}.) The impasse itself would be impossible}(Derrida 1993, p.21). An aporia cannot be located in a traditional Euclidian space due to a \say{lack of topographical conditions,} or, more radically, even a \say{lack of the topological condition itself}(Derrida 1993, p.21). This comment again marks what Derrida thinks to be the \say{limit} of mathematics, because it fails to locate the limit of time -- or, in this case, the limit between life and death— in either a Euclidian space or a topological space. 

\section{Arithmetic: The Number as \emph{Gramme}}

How do we perceive a \say{limit} that is both shapeless and impossible to be located in our real, three-dimensional space? The mathematics that belongs to the traditional system of metaphysics remains within the realm of points, lines, and planes, and is therefore unable to solve the problem. In this sense, mathematics has \say{limit} because its founding principles restrains itself from solving certain problems that are beyond its scope. But, as I argue, it is precisely the fact that mathematics has \say{limit} that offers an alternative way of thinking about the limit of time. The rest of this paper will primarily look at the three types of limit in three different branches of mathematics: arithmetic, calculus, and topology (with the aid of category theory, which could be seen a sub-field of topology).\footnote{Category theory, when it was first introduced in the 1940s, was created to solve problems in algebraic topology. As a general mathematical theory of structures, category theory has become an alternative to set theory as a foundation for mathematics. In this paper, I will mainly use category theory as an alternative way to think about topological spaces.} Let us briefly return to \say{Ousia and Gramme}, where Derrida suggests that seeing time in terms of \say{nows} represents a movement, a potentiality that is not arrested, an act that preserves the analogy of what he calls \say{gramme.} In arithmetic, as I will show below, the number \say{one} accurately mirrors the \say{now} as understood through the logic of \emph{gramme.} 

In the history of mathematics, numbers have been defined in many ways. One of the most ancient definition in Western thought was given by Euclid (1956: 277) in his \emph{Elements}: a number is a multitude consisting of units (Euclid 1956, p.277). This definition had been widely accepted until the modern period, when the 19th century logician Gottlob Frege proposed a new definition of number in terms of \say{equality} [\emph{Gleichzahligkeit}]. As he argues, the definition of numbers as a set or multitude is unable to cover 0 and 1. \say{The number which belongs to the concept F is the extension of the concept \say{equal to the concept}} (Frege 1980, p.85). For example, the number of the concept \say{author of \emph{Principe Mathematica}} is the extension of all concepts that are equinumerous to that concept. The statement \say{there are two authors of \emph{Principia Mathematica}} suggests that two objects fall in this class of concept $F$. The number 0, which implies that nothing falls under $F$, is understood as the number of the concept not being self-identical(Frege 1980, p.88).The problem of the Euclidian definition of number, according to Frege, is that it assumes a unity [\emph{Einheit}] under which numbers must be brought together, before defining numbers right away as a set of things. In other words, To think of \say{1} as a unit of counting is to take \say{one man} in the same way as \say{wise man,} that is, to think of \say{1} as an adjective that names a property of the thing it describes. Although Frege is not suggesting that there is no connection between \say{1} and \say{unity,} it is necessary to address the following problem before attempting to define numbers: are units identical to one another? 

Quoting Descartes, Frege argues that the plurality of number in things arises from their diversity.(Frege 1980, p.46). He agrees with W. S. Jevons that number is another name for diversity: unity is exact identity, and plurality arises from difference.(Frege 1980, p.46). He continues to quote the following example from Jevons, which illustrates both the identity and difference of numbers:

\begin{quote}
Whenever I use the symbol 5 I really mean

$1 + 1 + 1 + 1 + 1$,

and it is perfectly understood that each of these units is distinct from each other. if requisite I must mark them thus 

$1+1''+1'''+1''''+1'''''$ (Frege 1980, p.47).
\end{quote}

Are the \say{1}s really the same or are they intrinsically different? If they are indeed different, as Frege suggests, we would have to give up the equation of 1 = 1, and we should never be able to mark the same thing twice. \say{We must have identity – hence the 1; but we must have difference – hence the strokes; only unfortunately, the latter undo the work of the former.}(Frege 1980, p.48). As we can see, the same problem that we have encountered in the definition of time also arises here. The repetition of countless \say{nows,} which are understood as the \say{units} of time, also contains the difference that Frege mentions. On the one hand, every now is the same as the \say{now-no-longer} and the \say{now-not-yet}; on the other hand, it has its own radical singularity that is distinct from each other. Furthermore, the Aristotelean concept of time as a \say{number of change} that presides over the production of differences makes the counting of numbers possible. The \say{now} is simultaneously that which constitutes time (unity), and that which measures and divides time (difference). Attempting to overcome this difficulty, Frege suggests calling for assistance on time and space, because one point of space, or one \say{now} in time, is indistinguishable from another in a plane or a temporal line, except when conjoined as elements in a single total intuition. However, he quickly notices that number does not result from the mere division of the continuum of time and space: \say{time is only a psychological necessity for numbering, it has nothing to do with the concept of number.} Finally, \say{abstract number, then, is the \emph{empty form of difference}}(Frege 1980, p.56). For Frege, the concept of number has not yet been fixed, \say{but is only due to be determined in the light of our definition of numerical identity}(Frege 1980, p.74).

A detailed analysis of Frege’s solution to this problem and definition of number is less important for our purpose. The purpose is only to show that the evolution of the mathematical concept of \say{one} encounters a similar problem as that of \say{now}: as the foundational \say{unit} which constitutes and divides, both \say{one} and \say{now} face the challenge of being simultaneously same and different, universal and singular, present and absent. Aristotle defines time as \say{a number of change} [\emph{arithmos kineseos}] in respect to the earlier and later (Aristotle 1983, 219b1-2). The number 2, for example, is only a place holder for the integer between 1 and 3; similarly, every \say{now} is only relevant in its relationship to the \say{now} before [\emph{proteron}] and after [\emph{husteron}]. Time as \say{number of change} suggests that it presides over the production of differences, which makes the counting of numbers possible. Therefore, Aristotle’s \say{now} is not only something that constitutes time, but more importantly, that which measures and divides time (Aristotle 1983, 219b12). \say{The now is a link [syn{\'e}cheia] of time...for it links together past and future time, and is a limit of time, since it is a beginning of one and an end of another} (Aristotle 1983, 222a10-12). The limit of time, the \say{now,} connects time by dividing it, and vice versa. In a sense, the \say{now} functions as a punctuation which connects elements in a sentence by setting them apart, differentiating between past, present, and future phrases. 

In the section \say{Gramme and Number,} Derrida points out that, although time comes under the rubric of mathematics or arithmetic, it is not in nature a mathematical being (Derrida 1982, p.59). Derrida's conception of mathematics or arithmetic, I would argue, is Aristotelian; he might not notice the much stronger connection between \say{now} and the number \say{one}, as we have seen above. Aristotle think of number in two ways: numbering number and numbered number. \say{Numbering number} can be regarded as a universal form, in the sense that the number of a hundred horses and a hundred men is the same, whereas \say{numbered number} is a form embodied in matter, such as horses or men. Time, Aristotle claims, is a \say{numbered number} (Derrida 1982, p.59). Therefore, to think of time \say{arithmetically} \say{does not give time} because \say{[time] is as foreign to number itself, to the numbering number, as horses and men are different from the numbers that count them, and different from each other} (Derrida 1982, p.59). 

But what if each \say{numbering number} is, like each and every \say{now}, different from each other? Does there really exist a \say{limit} between \say{numbering number} and \say{numbered number}? If not, then there is no need to reject \emph{gramme} as a series of numbers, each of which would be \say{an arrested limit}(Derrida 1982, p.59). These numbers as limits are not \say{arrested,} closed points, but \say{empty forms of difference} through which similarity could be thought of. The limit of arithmetic is no longer a \say{multiplicity of points which are both origin and limit, beginning and end}(Derrida 1982, p.59), but possibly a \emph{gramme}, a potentiality, a movement, or something that might capture what Derrida calls a \say{trace.} Numbers share the same structure with \emph{gramme} as a potentiality produced by difference. Therefore, \say{[distinguishing] \emph{gramme} in general and the mathematical line,}(Derrida 1982, p.59, Aristotle 1983, 222a). seeing the mathematical line as a series of points, and thereby discounting the possibility understanding time in the mathematical nature, is assigning metaphysical thinking to mathematics, an act that illustrates the very tendency that Derrida aims at deconstructing.

\section{Calculus: The Self-Erasing Trace}

Derrida’s \say{trace} is a difficult, paradoxical concept, because it is impossible: that which enables something’s existence simultaneously annuls it. Every attempt to capture the trace (which simultaneously does and does not exist) through natural language is necessarily a failed attempt. Yet the trace is not unthinkable: in fact, mathematics teaches us that when something is unobservable, we could still \say{predict} it by getting infinitely close to it. How do we get infinitely close to something without actually reaching it? Here, we need to introduce the second type of \say{limit,} which is the limit in calculus. This concept defines the value that a function or sequence ‘approaches’ as the input approaches a certain value. Traditionally, an infinitesimal quantity is one which is smaller than any finite quantity without coinciding with zero, whereas an infinitesimal magnitude is what remains after a continuum has been subjected to an exhaustive analysis. 

Gottfried Leibniz was greatly preoccupied with the problem of the continuum, calling it one of the two \say{famous labyrinths where our reason very often goes astray}(the other being human freedom) (Leibniz 1951, p.53). Believing that he has found a solution to this problem, Leibniz writes in the preface to his \emph{New Essays}, \say{nothing takes place suddenly, and it is one of my great and best confirmed maxims that \emph{nature never makes leaps}: which is what I called the Law of Continuity} (Leibniz 1996, p.53). Leibniz‘s Law of Continuity claims that \say{when two instances or data approaches each other continuously, so that one at least passes into the other, it is necessary for their consequences or results (or the unknowns) to do so also} (Leibniz 1970, p.351). However, if a continuous entity has no gaps, the process of dividing it into ever smaller parts will never terminate in an indivisible point. Since real entities can only be either a simple unity or a multiplicity, and a continuum is not built from indivisible elements, Leibniz concludes that continua are not real \emph{entities} at all. 

In a sense, the limit in mathematics is also \say{impossible.} Suppose \emph{f} is a real-valued function and \emph{c} is a real number, the expression $\lim_{x \to c} f(x) = L$ means that $f(x)$ can be made to be as close to \emph{L} as possible by making \emph{x} sufficiently close to \emph{c}. However, by definition, \emph{x} can only infinitely approach \emph{c} without necessarily having to reach \emph{c}; consequently, $f(x)$ does not properly reach \emph{L}. Usually, limit helps us predict an answer when we cannot directly observe a function at a value. Therefore, with limit, calculus can attempt to answer seemingly impossible questions, such as when can a rectangle approximate a curve. However, limit only helps us predict an answer; it does not necessarily generate a function output at the very limit of the function. Consider the following example: 
\[
f(x)=(2x+1)(x-2)/(x-2)
\]
The limit of the function exists when $x=2$, because in that case the bottom will be 0, but we cannot divide any number by 0. Although we can predict that $f(2)=5$, $f(2)$ is never equal to 5, because $f(2)$ does not \say{exist}: when $x=2$, $f(x)$ would be divided by 0, the result of which would be meaningless, or it cannot be sensibly defined with real numbers and integers. The limit of this function tells us that we can only get as close as possible to the point of $f(2)$, but as soon as we arrive there, the point erases itself, leaving a gap in the graph of the function. The limit in calculus is, like the limit of time, only produced by its own erasure. 

Derrida sees time as \say{the name of the limits within which the \emph{gramme} is thus comprehended, and, along with the \emph{gramme}, the possibility of the trace in general} (Derrida 1982, p.60). Time is unable to arrest each concrete moment that has a beginning and an end. It is only able to catch its own trace, just as the form of calculus – in which infinitesimals are represented in terms of limit procedures – makes it possible to grasp the output of a function at its limit. But \say{the trace is produced as its own erasure}(Derrida 1982, p.65); it is neither perceptible nor imperceptible. We can never perceive the trace that time is able to catch as a concrete entity, but only predict its existence. At the limit of time, trace can only exist as something else. 

The second form of mathematical limit provides a rigorous mathematical model to formally think about Derrida's \say{trace,} a paradoxical concept that cannot be captured by natural language. If a \say{limit} in our natural language functions as a line, a threshold that separates two individual concepts, a \say{limit} in calculus functions as an arrow: it points at a certain unattainable answer, it offers the best solution possible without providing a final answer. But in our attempt to locate the \say{limit of time,} we encounter yet another problem: can we obtain a best approximate answer to the question, \say{what is time?} Is there a satisfying answer to the metaphysical essence of time, even without knowing what happens when I reach the limit of my time, what happens after my time? The shapeless, timeless concept of death is indeed horrifying: is it itself a limit, a border? Or, does it even have a border? 

\section{Topology: The ``Shape" of Aporia}

25 years after \say{Ousia and Gramme,} in \emph{Aporias}, Derrida discusses the problems of the limit again, focusing on the central field of non-limitation. The word \say{aporia} that appears in \emph{Physics IV} reconstitutes the aporia of time via an exoteric logic (Derrida 1993, p.13, Aristotle 1843, 217b). In this much later essay, Derrida suggests that he was trying to demonstrate that the philosophical tradition from Kant to Hegel has only inherited the Aristotelean aporetic, namely, that now is and is not what it is (it only \say{scarcely} is what it is)(Derrida 1982, p.14). If \say{Ousia and Gramme} confirms the hegemony of the \say{vulgar} concept of time insofar as it privileges the \say{now} in the philosophical tradition, \emph{Aporias} further questions: \say{what if there was no other concept of time than the one that Heidegger calls \say{vulgar}? What if, consequently, opposing another concept to the \say{vulgar} concept were itself impracticable, nonviable, and impossible? What if it was the same for death, for a vulgar concept of death?}(Derrida 1982, p.14) Both questions about time and death have been summarized as a single question of the aporia, the impassible passage that marks the limit of time. If in \say{Ousia and Gramme,} the limit of time is marked by the trace of every \say{now}, then in \emph{Aporias}, it becomes a question between life and death. The task of \say{Ousia and Gramme} is to present a possible shape and location of the limit of time, whereas in \emph{Aporias}, the task becomes to transgress such limit (or, rather, to wait for it). 

In \emph{Being and Time}, Heidegger defines death as \say{the possibility of the pure and simple impossibility for Dasein} (Heidegger 2008, p.250). Indeed, death is an \say{occurrence} of this possible impossibility, because at the very the moment that \emph{Dasein} \say{dies,} it is no longer possible [\emph{Daseinunm{\"o}glichkeit}]. The paradoxical structure of aporia -- its inaccessible passage -- stems from the fact that there is no limit to cross: there is no opposition between two sides, or, the line is too porous, too permeable, and impossible to locate (Derrida 1982, p.20). If so, can we still reduce death to the crossing of a border? Where do we find this unique aporia, or, in other words, where does one \say{expect death,} as expecting the impossibility of the possible?

Derrida's analysis of the aporia consists in three parts. He defines the three types of border limits: first, anthropologico-cultural borders that separate territories, nations, cultures, and languages; second, those that separate domains of discourses such as philosophies, knowledges, disciplines of research, or, as Derrida later refers to, \say{the delimitations of the problematic closure}(Derrida 1982, p.23); third, the lines of separation, demarcation, or opposition between concepts or terms – this is what Derrida calls \say{the conceptual demarcation} (Derrida 1982, p.73). Then he attempts to situate the possible impossibility of \say{my death} among these marked borders and limits. The limit between life and death, according to Derrida, is not an \say{absolutely pure and rigorously uncrossable limit}(Derrida 1982, p.75), but the possibility of an impossible passage. The problem, as Derrida sees it, is the following: if death is the ultimate limit, then how can we place the edge of a limited life against a potentially illimitable death to even form a threshold? To think about \say{the other side} after death is only possible if one \say{has elaborated a concept of the ontological essence of death and if one remembers the possibility of being of every \emph{Dasein} is engaged, invested, and inscribed in the phenomenon of death}(Derrida 1982, p.53).In other words, the thinking of \say{the other side} can only start from \say{this side.} Furthermore, the metaphysics of death cannot claim to have any coherence or rigorous specificity, because we only have access to the process of dying, not death itself. 

However, Derrida's conclusion does not tell us where we should await our own death: how do we experience the uncrossable limit? For Derrida, there is not even any \emph{space} for the aporia, \say{because of a lack of the topographical condition or, more radically, because of a lack of the topological condition itself}(Derrida 1982, p.21). Interestingly, he is using the two concepts, topography and topology, interchangeably here. Later, he even suggests that topology is the pre-condition of topography. To a certain degree, \say{aporias} suffers from this conceptual confusion between topography and topology: while the former relies on a geometric, topographic vocabulary, the latter, especially in the mathematical sense, is fundamentally abstract and shapeless.  In William Watkin’s \say{Derrida's Limits -- Aporias between \say{Ousia and Gramme},} he points out that the habitual reliance on a geometric vocabulary in metaphysics is one of the main causes of the aporia of the limit that Derrida speaks about (Watkin 2010, p.124). Watkin’s argument is that, while the aporia between life and death cannot be crossed, because the two are anti-morphic (that is, they do not share the same shape), the study of topology makes the border between life and death crossable, since both sides share certain common characteristics. 

Watkin aims at tackling two paradoxical questions raised by Derrida, namely, \say{can death be reduced to some line crossing, to a departure, to a separation, to a step, and therefore a decease?} and its interro-denegative, \say{is not death, like decease, the crossing of a border, that is, a trespassing on death, an overstepping or a transgression?} (Derrida 1983, p.6) There exists \say{something} between life and death that separates the two, but the possibility of drawing a limit between life and death is impossible because death is \say{illimitable.} Therefore, to cross from life to death simultaneously is and is not the crossing of a border or a limit, because the limit is and is not crossable. Watkin, following the lead of Derrida, continues to pursue the question of the uncrossable limit between different fields, subjects, disciplines, et cetera. He argues that the limit \say{both is and is not aporia} because the limit depends on its own logic: it must open itself to the field of non-field in order to define and form the field (Watkin 2010, p.117). The limit \say{touches} both that which is limited and that which is illimitable, the latter of which it cannot cross. Therefore, Watkin suggests that the limit should be defined as the \say{possibility of the impossibility of being’s continuation and pure extension in all directions including extension into itself} (Watkin 2010, p.117).

Then, Watkin associates geometric topographic limit with closed self-identical fields whose borders touch each other without remainder or gaps, and topological limit with convergences, connections, and continuities. Since we overly rely on the geometric topographic vocabulary as we speak of philosophical concepts, such as the idea of life and death, we tend to assume such concepts as \say{closed self-identical fields whose borders touch each other by something we conceive of as a kind of line which touches equally both sides of the divide} (Watkin 2010, p.124). Since mortal life is limited, while death is potentially illimitable, the limit between the two cannot be crossed because life and death are \say{anti-morphic.} Inspired by the famous example of the mug-doughnut in the study of topology, in which a mug is the equivalent of a doughnut because of their shared commonality across basic classes of homotopy (namely, holes and tails), Watkin suggests that the aporia between life and death should be thought of as a topological border, namely a continuum shared by the two concepts in terms of their commonalities: mortality, temporality, and so on. He concludes with a witty analogy: \say{life as limit bears within its death illimitability which produces a finitude to life which death as finitude of life carries inside itself in the form of life as the perpetually open. Life is the same class as death in that it has one hole and no tail. Life is death-shaped; death is a torus. Death is a doughnut}(Watkin 2010, p.126).

Watkin's example of topological homeomorphism is, of course, correct; however, it is rather arbitrary to think of life and death as two \say{homotopy equivalent} in topological spaces. The section entitled \say{Topology: The Shape of Death} presupposes that death has a shape, either as the traditional topographical markings of the limit as a line or a circle, or as the torus-doughnut that Watkin depicts. Similarly, Derrida see death as taking a figure (Derrida 1993, p.58). It has a \say{privileged form, the crossing of a line} between existence and non-existence, \emph{Dasein} and non-\emph{Dasein}. At the very end of \say{Ousia and Gramme,} Derrida illustrates the potential shape of the trace of time: \say{a writing exceeding everything that the history of metaphysics has comprehended in the form of the Aristotelian \emph{gramme}, in its point, in its line, in its circle, in its time, and in its space} (Derrida 1982, p.58). But do we know really so much about time or its limit that we can depict its form, draw its border, and mark it as it crosses the line? Quoting Heidegger, Derrida stresses the concept of \say{awaiting oneself} [\emph{s'at-tendre}] in the existential analysis of Dasein and its being-towards death: \say{death is a possibility-of-being that Dasein itself has to take over [zu {\"u}bernehmen] in every case. With death, Dasein awaits itself [\emph{s'at-tend lui-m{\^e}me, steht sich… bevor}] in its own most potentiality-for-being} (Heidegger 2008, p.250, Derrida 1982, p.64). Here, the \say{awaiting} seems to be impossible without spatial relationships: one await oneself \say{in} oneself, or \say{directed towards} some other; with death, Dasein awaits \say{in front of} itself, \say{before} itself like in front of a mirror (Derrida 1982, pp.65-66).\footnote{In \emph{Being and Time}, Heidegger has already considered the spatial character of prepositions, which suggest that, in time, there is always already a certain spatialization that hides the matter at hand. Here, Derrida is working along the same lines.} If, as Watkin says, \say{the shape of life and the shape of death are homotopic in some areas}(Watkin 2010, p.126), then the two must share several homeomorphism classes. But are life and death even shapes, to begin with? What kind of spatial or temporal structures do they have, and if so, is it possible to trace their border? How do we know that death is really just a doughnut?

\section{Category Theory: the Shapeless Shape}

The problem of thinking life and death as some specific topographical shapes is that it is not general enough. Instead, I suggest that we imagine this limit between life and death -- two \say{categories} of which we still know too little to define -- as a limit in category theory, an established field of mathematics. What Watkin's example offers can be translated as \say{plugging} in one object (whose shape we already know -- a torus) into one of our categories and generating a certain answer, which is confined \say{locally.} In order to study the two categories, we still need to take one step further and reach the same level of abstraction. Assuming a topological \say{shape} of death is like studying the ontology of time when we already have a (Kantian) \emph{a priori} intuition of time. Such practice would only be considered tautological. Therefore, to simply \say{glue} the two sides it separates back together, as what Watkin suggests in his \say{double torus} example, is not sufficient for this problem. To tackle the problem of the \say{impossible possibility} of Derrida's limit, we need to use the language of general topological space that can be applied to categories even when we do not know the \say{shape} of them. I argue that the \say{limit} between life and death is a \say{shapeless shape} that, although oftentimes invisible and unimaginable, nevertheless tells us more about topological spaces that we have no prior knowledge about. To understand this \say{shapeless shape,} we first need to know why this limit is \say{shapeless} -- that is, it is neither a line, nor a hole. 

The study of topology is concerned with the properties of spaces that are preserved under continuous deformations, such as the act of stretching a \say{mug} into a \say{doughnut.} Unlike the study of geometry or calculus, topology requires very little structure: all it needs are two notions, \say{points} and \say{nearby points.} The concept of topology is indeed helpful because we no longer need to think of concepts in terms of their spatial structures; all we need are some points in a so-called \say{open set.} Therefore, if we follow Watkin’s lead and think of life and death as two topological concepts, we should stop imagining them as \say{shapes,} as we know too little about life and death to assign them geometric topographic structures. What we know about these concepts is limited by fundamental logical relationship such as union, intersection, complement, and differences. We can also look at them from local perspectives, such as that of a point. Since the study of topology helps us gain a better understanding of abstract concepts when we cannot even begin to imagine its shape, we may have to stop thinking about the \say{shape} of death, as Watkin’s chapter entitles, but indeed, the limit of it. 

In topology, the idea of \say{limit} helps establishing a sort of relationship between points. Let us first look at the definition of \say{limit} in a topological space. The limit of a sequence of points (indexed by the natural numbers)$x_1, x_2, ..., x_n,...$ is a point $x$ such that:

\begin{quote}
- With any open set $U$ which contains $x$:

- There is a number $n$ such that for any number $N$ greater than $n$,

- $X_N$ is in $U$.
\end{quote}

In a \say{vulgar} sense, we can learn from this definition that in a topological space, there is always another point that has some relationship with the given point. Instead of a set of scattered points floating around the space, the limit of a topological space gathers points, bringing them together. At this point, the \say{limit} neither demarcates the border of a closed self-identical field, nor exists as a continuum between two spaces, since there are no longer shapes, fields, or spaces. In topology, \say{limit} depicts the arbitrary \say{closeness} between two points without measuring the distance between them. It does not even have a geometric structure: It is neither a line nor an arrow. In fact, we know so little about the limit that we cannot even draw a picture of it. Then what is a limit? Mathematicians use limit to establish connections between existing points and re-construct a topological space with it. In a sense, we can say that the limit functions as an invisible \say{link} that connects, twists, folds, creates, and constructs. 

If this \say{shapeless} topological limit is too abstract to be translated into philosophical language that helps us cross the uncrossable limit between life and death, category theory as a mathematical tool helps us delineate this shapeless shape. As we have seen earlier, in mathematics, a limit is usually understood as the value that a function \say{approaches} as the input approaches some value. This idea may be applied to functions, sequences, or geometric structures. But what if we want to look at more general topological spaces beyond manifolds that can be embedded (translated, drawn out) in three-dimensional real spaces? In this case, the traditional numerical invariants used to study these spaces would be rendered meaningless, and it therefore becomes necessary to study these spaces taken as a whole (Volkert 2002, p.284). During the early 1940s, Samuel Eilenberg and Saunders Mac Lane collaborated to look for a tool to study such these groups of spaces (Eilenberg and Mac Lane 1942, p.759). The two mathematicians systematically expose the concept of ‘category,’ a term borrowed from the philosophy of Aristotle and Kant (Eilenberg and Mac Lane 1971, 29f.). Here, I will briefly introduce the concept of category with the help of the graph below. By definition, a category is a labeled direct graph with nodes (also called objects) and arrows (formally morphisms). As shown in the graph below, $X, Y$ and $Z$  are objects, whereas $f,g$ and $g \circ f$  are morphisms. The collection of $X, Y, Z, f, g,$ $g\circ f$ such forms a category.

\[
\begin{tikzcd}
X \ar{r}{g\circ f}\ar[swap]{d}{f} & Y \\
Z \ar[swap]{ur}{g} & 
\end{tikzcd}
\]

With categories, mathematicians are now able to axiomatically capture common features of groups of mathematical structures. By studying the structure-preserving mappings between objects or groups of objects, one could \say{categorize} many areas of mathematical study such as sets (a collection of distinct objects), groups (algebraic structures that include a set and an operation, such as addition or subtraction), and topologies (a set of points along with a set of neighborhoods for each points). To \say{categorize} is to \say{assemble} groups of objects into a category, displaying the relationships between these different objects, as opposed to considering them as \say{disembodied} pieces. Although each category has distinct properties, in the definition of a category, objects are considered atomic, therefore we do not know whether an object is a set, a topology, or something else. Since category theory studies spatial objects without referring to their internal structure, it can provide a higher level of abstraction, making it possible to find \say{universal properties} that determine these objects of interest. To find universal properties is another way of saying that, instead of considering the \say{internal structure} of an object, we will instead look at the object as a category, defined in terms of its spatial relationship with its neighbors. If we could apply categorical thinking to the philosophy of limit, we would be able to find another way to reflect on the unknown category (death) through the already-known category (life), without necessarily \say{crossing} the uncrossable border between the two. 

The categorical notion of limit, introduced in order to describe objects in a purely categorical way, captures the essential properties of universal constructions in category theory. Limit, as a universal property in category theory, is defined differently as in calculus or arithmetic, not because they are intrinsically two concepts that have nothing in common,\footnote{Even though the notion of limits in category theory was indeed inspired by the classical notion of a limit.} but because at a higher abstraction level, the definition of a limit can no longer rely on specific features of an object, such as its shape or location. Even so, the limit in category theory still has a \say{shape,} although not a geometric one: it is a \say{cone.} This \say{shape} captures the interrelationship between one object and other objects. For example, let us say we have two categories $C$ and $J$ , and $F$ is a functor from $J$ to $C$. (That is to say, $C$ is a directed graph, and $J$ carves out a subgraph of $C$.)
\[
F:J\rightarrow C
\]
The function of morphisms in $F$ maps every object in category $C$ to an object in category $J$, thus establishing a relationship between the two categories. In mathematical language, everything in $C$ \say{commutes} to everything in $J$. Then in order to draw the shape of $F$, we have to \say{map} every object in category $C$ to category $J$, which resembles the process of \say{carving out} the shape (that is, a collection of nodes and arrows) of $J$ in $C$. With this information, we can draw a diagram indexing all the objects and morphisms in $C$ patterned on $J$. Then we pick some random object $c$ in the category $C$ as the apex of our \say{cone.} The groups of morphisms connecting $c$ to every functor from $J$ to $C$ shapes our \say{cone.} Despite its name, it is not always possible to locate, measure, and draw out its distinct shape in a Cartesian space, because we do not necessarily know the shape of our objects. We could, however, provide a concrete example in the Cartesian space: if $J$ is a circle, and $C$ is the line or interval between $0$ and $1$. Then the \say{cone} in this example actually has the shape of a cone, with a flat, circular base that tapers smoothly to its vertex. 

Applying topological thinking to the study of philosophy could greatly broaden our scope: we no longer need to know if death is a doughnut, or if it has a shape at all. According to Derrida, to answer the question \say{is my death possible,} we first need to know how to \say{cross} from life to death (Derrida 1993, p.21). Yet each death is so absolutely singular and \say{unassignable} that it necessarily remains unknown (Derrida 1993, p.22). It is impossible to demarcate a \say{limit} in the sense of \say{border} or \emph{gramme} between life and death, because we know neither the \say{shape} nor \say{location} of the latter—we don’t even know if it is an \say{object.} Therefore, Derrida calls that which lies between life and death \say{aporia,} because it is an impossible yet necessary passage that cannot be crossed. Topological (and categorical) thinking solves this problem by creating another \say{limit} that essentially draws the map of our \say{commute} to the unknown. We are no longer required to \say{cross the limit} to reach the other side; all we need to know is how to find a different \say{limit} that depicts the relationship between life (the object that we already know) and death (the unknown \say{category}). By doing so, we will be able to learn something about the shapeless, traceless, limitless \say{other} that lies beyond the uncrossable aporia.

\end{document}